\input amstex
\documentstyle{amsppt}
%
%
%
%
%
%
%\NoBlackBoxes
\NoRunningHeads
\TagsOnRight
\magnification=1200
%\loadbold
%
\font\tensmc=cmcsc10
%
%%%%%%%%%%%%%%%%%%%%%%%%%%%%%%%%%%%%%%%%%%%%%%%%%%%%%%%%%%%%%%%%%%%%%%%%%%%%
%%%%%%%%%% TeX macros.
\define\thismonth{\ifcase\month % case 0 --- impossible!
  \or January\or February\or March\or April\or May\or June%
  \or July\or August\or September\or October\or November%
  \or December\fi}
\newcount\formnumb
\formnumb = 0
\def\calcformnumb#1{{\count255=\formnumb\advance\count255 by #1
\the\count255}}
\def\nfn{{\global\advance\formnumb by 1}\the\formnumb}
\def\pfn#1{{\bf (\calcformnumb{#1})}}
\def\romanitem#1{\itemitem{{\rm #1}}}
\def\lsc{l\.s\.c\.\ }
\def\norm{\|\cdot\|}
\def\dist{\roman{dist}}
\def\diam{\roman{diam}}
\def\coenv{\roman{co}\,}
\def\arginf#1{{{\frak{S}}_{{\frak{Inf}}}({#1})}}
\def\unifder{{\Cal {C}}^{1,u}}
\def\bdunifder{{\Cal {C}}_{\Cal {B}}^{1,u}}
\def\holder#1{{\Cal {C}}^{1,#1}}
\def\bdholder#1{{\Cal {C}}_{\Cal {B}}^{1,#1}}
\def\extR{\Bbb{R}\cup\{+\infty\}}
\def\infconv#1#2{\big( #1\phantom{[} \square\phantom{]} #2 \big) }
\def\Deltconv#1#2{#1\blacktriangle #2}
\def\cqd{\enspace\lower.5ex\hbox{$\square$}}
%
%%%%%%%%%%%%%%%%%%%%%%%%%%%%%%%%%%%%%%%%%%%%%%%%%%%%%%%%%%%%%%%%%%%%%%%%%%%%%
%%%%%%%%%% Title, abstract,...
\topmatter

\title 
On regularization in superreflexive Banach spaces 
by infimal convolution formulas
%\\
%
%
%$\roman{(Preliminary\ version)}$
\endtitle

\author 
Manuel Cepedello Boiso
\endauthor

\address
Equipe d'Analyse, Universit\'e Pierre et Marie Curie--Paris 6, Paris.
\newline
\indent
Departamento de  An\'alisis Matem\'atico, Universidad de Sevilla, Sevilla.
\endaddress
\curraddr
Department of Mathematics, University of Missouri-Columbia, Columbia.
\endcurraddr

\email
manuel\@lebesgue.math.missouri.edu, cepedel\@ccr.jussieu.fr,  
boiso\@cica.es
\endemail

%  Math Subject Classifications 

\abstract 
We present here a new method for approximating functions defined on 
superreflexive Banach spaces by differentiable functions 
with $\alpha$-H\"older derivatives (for some $0<\alpha\leq 1$). 
The smooth approximation is given by means of an explicit formula
enjoying good properties from the minimization point of view.
For instance, for any function $f$ which is bounded below and
uniformly continuous on bounded sets this formula gives 
a sequence of $\Delta$-convex 
${\Cal{C}}^{1,\alpha}$ functions
converging uniformly on bounded sets to $f$
and preserving the infimum and the set of minimizers of $f$.
The techniques we develop are based on the use of 
{\sl extended inf-convolution} formulas and convexity
properties such as the preservation of smoothness for the convex envelope of 
certain differentiable functions. 
\endabstract

\date
%
%{\thismonth} {\number\day}, {\number\year}
%
May 28, 1997
\enddate


\thanks
The author was supported by a FPU Grant of the Spanish {\it Ministerio de 
Educaci\'on y Ciencia}.   
\endthanks

\endtopmatter
%%%%%%%%%%%%%%%%%%%%%%%%%%%%%%%%%%%%%%%%%%%%%%%%%%%%%%%%%%%%%%%%%%%%%%%%%%%%
%%%%%%%% Text
\document
%
%
%
%%%%%%%%
\head 0.
Introduction and Preliminaries
\endhead
%%%%%%%%
\par
This paper introduces an explicit regularization procedure for functions 
defined on superreflexive Banach spaces. For any bounded below \lsc 
(resp. uniformly continuous on bounded sets) function $f$ on a superreflexive 
Banach space $X$ we give by means of a ``standard'' formula a sequence of 
${\Cal{C}}^{1,\alpha}$-smooth functions
converging pointwise (resp. uniformly
on bounded sets) to $f$ (where $0<\alpha\leq 1$ only depends on $X$).
Under some additional conditions, the convergence of the sequence of 
approximate functions
is uniform on the whole space $X$. Moreover, the approximate functions 
preserve
the infimum and the set of minimizers of $f$. We remark that these features 
altogether cannot be easily 
obtained from regularization methods like the 
{\it smooth partitions of the unity} techniques
(for a detailed study of this topic we refer to Chapter VIII.3 of \cite{DGZ},
the references therein and \cite{Fr}) or other results that only ensure the
existence of smooth approximates (for instance, see \cite{DFH}).
\par
In Hilbert spaces, our work is closely linked with the 
{\it Lasry-Lions approximation method}  (introduced in \cite{LL} 
and subsequently 
studied by several authors, such as \cite{AA}) and its more general version 
given by T. Str\"omberg in \cite{St$_2$}. Actually, we improve the results
of \cite{St$_2$} in the superreflexive case by providing the best 
uniformly smooth approximation possible for this setting. Nonetheless, 
we want to remark that the approximate functions explained 
herein cannot be reduced to
those of Str\"omberg (or Lasry-Lions approximates in Hilbert spaces); 
we refer to the remark after {\bf Proposition 8} for a more precise 
explanation. 
Our approach for smooth regularization in non-Hilbert spaces comes from 
two main facts:
the density of the linear span of the convex functions (studied in \cite{C})
and the smoothness of the convex envelope of a ``somehow'' smooth function.
In this direction, we also present more general versions of certain results in
\cite{GR} for infinite dimensional Banach spaces.
\par
This paper is organized in the following way.
Our main result of this paper, {\bf Theorem 1}, and several corollaries 
are explained in
{\smc Section 1}. The proof of {\bf Theorem 1} is showed in {\smc Section 4}
with the tools provided by sections 2 and 3. {\smc Section 2} deals
with the existence of approximates for a given function $f$ 
using some results on extended inf-convolution formulas.
{\smc Section 3}
develops a procedure 
for regularizating certain $\Delta$-convex approximates.
This procedure is based on the smoothness of the convex envelope of 
certain ``somehow'' smooth functions.
\par
%
%%%%%%%%
\subhead\nofrills
{\smc Notation: }
\endsubhead
%%%%%%%%
%
In what follows, $X$ denotes a Banach space and $\norm$ an equivalent 
norm on $X$. By $B_X$ we mean the unit closed ball of $X$ under 
the norm $\norm$ and by $B_X(r)$ the closed ball of radius $r>0$. 
A function $f:X\to \extR$ is called 
{\sl proper} if $f\not\equiv +\infty$ and  $\arginf{f}$ 
is the (possibly empty) set
$\{x\in X : f(x)=\inf f\}$.
We will deal with the pointwise, compact, uniform on
bounded sets and uniform on $X$ convergence in the set 
of lower semi-continuous
(in short, \lsc) functions on $X$, abbreviated respectively by 
$\tau_p$, $\tau_K$, $\tau_b$ and $\tau_u$.
\par
A function defined on $X$ is called {\sl $\Delta$-convex} if 
it can be expressed as the difference of two continuous convex functions.
The convex envelope $\coenv{f}$  of a function $f:X\to \extR$
is defined as the
greatest proper convex \lsc function below $f$ (if there exists 
a convex minorant of $f$). 
The explicit value of the convex envelope of $f$ at a point $x\in X$ 
is given by the formula 
$$
(\coenv{f})(x)\!=\!\inf_{n\in\Bbb N}\!\Bigg\{\!\sum_{i=1}^n \lambda_i f(x_i): 
x\!=\!\sum\limits_{i=1}^n\lambda_i x_i,\, 
\sum\limits_{i=1}^{n} \lambda_i\!=\!1,\,
\big(x_i,\lambda_i\big)_{i=1}^n\! \subset \! (X\times{\Bbb {R}}_{+})\! 
\Bigg\}.
\eqno(\nfn)
$$
\par
Unless stated otherwise, differentiability will be understood in the 
{\sl Fr\'echet sense}.
The following notation is used throughout this work. By $\unifder(X)$ 
(respectively $\bdunifder(X)$) we understand the set of 
differentiable functions defined on $X$ with
uniformly continuous (resp. uniformly continuous on bounded sets) derivative. 
Similarly,  
$\holder{\alpha}(X)$ (resp. $\bdholder{\alpha}(X)$) 
stands for the class of functions on $X$ having $\alpha$-H\"older continuous 
(resp. $\alpha$-H\"older continuous on bounded sets) 
derivative ($0<\alpha\leq 1$).
\par
%
%%%%%%%%
\head 1. 
The main result
\endhead
%%%%%%%%
\par
We begin by stating the main result of this work.
\proclaim{Theorem 1}
Let $p>1$, 
$X$ be a Banach space and $\norm$ be an equivalent norm on $X$ 
which is locally uniformly convex and uniformly smooth. 
For any proper lower semi-continuous bounded below function 
$f:X\to\extR$, consider 
the sequence of $\Delta$-convex functions given
by the formula
$$
\Delta_n^p{f}:= \coenv{g_n^p} -  2^{p-1}n\norm^p 
\quad
(n\in \Bbb N),
$$
where $g_n^p$ at a point $x\in X$ is defined as  
$$
g_n^p (x)
:= 
\inf_{y\in X} 
\Big\{ 
f(y)+2^{p-1}n\|x\|^p + 2^{p-1}n\|y\|^p - n\|x+y\|^p 
\Big\}
+ 
2^{p-1}n\|x\|^p.
$$
Then the following assertions are satisfied:
%
%\roster
%
\romanitem{(i)}  For all $n$, 
$\inf f \leq \Delta_n^p{f} \leq f$ and $\arginf{\Delta_n^p{f}} =\arginf{f}$.
\romanitem{(ii)} $(\Delta_n^p{f})_{n\in\Bbb N} \subset \bdunifder(X)$ and  
$(\Delta_n^p{f})_{n\in\Bbb N} \subset \bdholder{\alpha} (X)$ 
provided that the modulus of smoothness of
the norm $\norm$ is of power type $1+\alpha$; actually, we have that
$(\Delta_n^{1+\alpha}{f})_{n\in\Bbb N} \subset \holder{\alpha}(X)$.
\romanitem{(iii)} $\Delta_n^p{f} @>\tau_p>n> f$ pointwise and $
\Delta_n^p{f} @>\tau_K>n> f$  if $f:X\to\Bbb R$ is continuous.
\par
If moreover the norm $\norm$ is uniformly convex then 
\romanitem{(iv)} $\Delta_n^p{f} @>\tau_b>n> f$ whenever $f$ is uniformly 
continuous on bounded sets.
\romanitem{(v)} $\Delta_n^p{f} @>\tau_u>n> f$ provided that $f$ 
is uniformly continuous on $X$ (not necessarily bounded below) 
and the modulus of convexity of the norm $\norm$ is of power type $p$
($p\geq 2$).
%
%\endroster
%
\endproclaim
\par
\remark{Remark}
It is well-known that the existence of a uniformly smooth norm $\norm$
on a Banach space $X$ implies the superreflexivity of $X$ (and reciprocally, 
the articles  \cite{E} of P. Enflo and \cite{Pi} of G. Pisier 
tell us that any superreflexive Banach space 
admits an equivalent uniformly smooth norm). 
Similarly, we want to point out
that the conclusions of {\bf Theorem 1} cannot be expected outside the 
superreflexive setting. 
\par
First, the $\tau_b$-density of the set of $\Delta$-convex functions 
defined on $X$
in the set of functions on $X$ that are uniformly continuous on bounded sets 
is equivalent to the
superreflexivity of the Banach space $X$(as it was proved in \cite{C}).
On the other hand, the existence of ${\Cal {C}}^{1,\alpha}$ bump functions
(for some $0<\alpha\leq 1$) on $X$ implies the existence of 
an equivalent norm $\norm$
on $X$ with modulus of smoothness of power type $1+\alpha$ (see 
Theorem V.3.1. of \cite{DGZ}). 
\par
\endremark
\par
\remark{Remark}
The optimal application of {\bf Theorem 1} is achieved when we consider
a Hilbertian norm $\norm$. In this case, taking $p=2$ in {\bf Theorem 1}
we obtain similar
approximation results as those given by the Lasry-Lions approximation method
(see \cite{LL}). Nevertheless, the different sequences of approximates are not
the same even in this setting (see remark after {\bf Proposition 8}).
\endremark
We proceed to state some corollaries to {\bf Theorem 1}. They are related with
certain results known on a superreflexive Banach spaces from the existence of
{\sl smooth partitions of the unity} (see Theorem VIII.3.2 in \cite{DGZ}).
Their proof is easily obtained appealing to {\bf Theorem 1} and 
Pisier's renorming
Theorem (the original proof can be found in \cite{P}; we refer to \cite{L} 
for a simpler and more geometrical proof).
\par
The first corollary improves Corollary 1 of \cite{St$_2$} for superreflexive
Banach spaces.
\par
\proclaim{Corollary 2}
Let $X$ be a superreflexive Banach space. 
Then there exists some $0<\alpha\leq 1$ such that 
any non-empty closed set $F$ of $X$
is the set of zeros of a $\Delta$-convex $\Cal{C}^{1,\alpha}$-differentiable
 function on $X$. 
Moreover, $F$ is the limit for the Hausdorff distance of a sequence
of sets $\frak {S}_n=\{x\in X: f_n(x)<\sigma_n\in\Bbb R\}$ ($n\in\Bbb N$)
where 
the functions $(f_n)_n$ are $\Delta$-convex and 
in $\bdholder{\alpha}(X)$.
\endproclaim
\par
\demo{Proof of Corollary 2}
For a superreflexive Banach space $X$, Pisier's renorming Theorem 
ensures the existence of an equivalent norm $\norm$ on $X$ with modulus of 
smoothness of power type $q$ ($1< q\leq 2$). Given a closed set $F$ in $X$,
consider the proper function $d$ defined at a point $x\in X$ as 
$d(x):=\dist(x,F)=\inf_{y\in F}\|x-y\|$ ($d$ is proper because 
$F$ is not empty).
By {\bf Theorem 1}(i)--(ii) we have that the function
$\Delta_1^{q}{d}$ is $\Delta$-convex, $\Cal{C}^{1,q-1}$-differentiable and 
satisfies that 
$\arginf{\Delta_1^q(d)}=\arginf{d}=F$.
\par
Moreover, using Asplund averaging technique 
(see Proposition IV.5.2 of \cite{DGZ}),
we can assume that the
modulus of convexity of the norm $\norm$ is in addition 
of power type $p$ 
(for some $p\geq 2$). 
Since $d$ is Lipschitz continuous on $X$,
from {\bf Theorem 1}(iv) it follows for every $n$ that 
$F\subseteq\{\Delta_n^{p}{d}(x)<\frac {1}{n}:x\in X\}:=\frak{S}_n$, where 
$\Delta_n^p{d}$ is a $\Delta$-convex $\bdholder{q-1}$-differentiable 
function and $(\frak{S}_n)_n$ converges to $F$
for the Hausdorff distance.
\cqd
\enddemo
\par
The next corollary gives a slightly stronger version of some others 
approximation results obtained by using partition of the unity techniques
(for instance, see Theorem 1 of \cite{NS}).
\par
\proclaim{Corollary 3}
For any superreflexive Banach space $X$ there is $0<\alpha\leq 1$ so that
for every uniformly continuous on bounded sets (resp. uniformly continuous)
function on $X$ one has the following: 
$f$ is the uniform limit on any fixed bounded set $B$ of $X$ 
(resp. on $X$)
of a sequence of $\Delta$-convex
$\holder{\alpha}$-differentiable 
(resp. $\bdholder{\alpha}$-differentiable)
functions
having the same infimum and set of minimizers on $B$ as $f$. 
\endproclaim 
\par
\demo{Proof of the Corollary 3}
Appealing again to Pisier's renorming Theorem for superreflexive 
Banach spaces,
we can suppose that there is an equivalent norm $\norm$ on $X$ 
with modulus of
smoothness of power type $q$ ($1<q\leq 2$). Fix some bounded set $B$ of
$X$ and define $\tilde{f}:=\max\{f,\inf_B f\}$. 
Since $f$ is uniformly continuous on $B$, we have that $\inf_B f>-\infty$.
Therefore, $\tilde{f}$ is 
uniformly continuous on bounded sets and bounded below.
Note that trivially $\tilde{f}(x)=f(x)$ for all $x\in B$ and then
the infimum and set of minimizers on $B$ of $f$ and $\tilde{f}$ are the same.
Hence,
{\bf Theorem 1}(ii) and (vi) tell us that the sequence
$(\Delta_n^q \tilde{f})_n$ satisfies the required conditions of the claim
for $\alpha=q-1$. 
If $f$ is uniform continuous on $X$, the proof of Corollary 3 follows the same
lines, using the existence on  
$X$ of an equivalent norm $\norm$ with non-trivial moduli of convexity and
smoothness and {\bf Theorem 1}(v).
\cqd
\enddemo
\par
The last corollary is an extension of Remark (viii) in \cite{LL}. It deals
with the property of extending and regularizing functions defined
on subsets of superreflexive Banach spaces to the whole space.
\par
\proclaim{Corollary 4}
Let $X$ be a superreflexive Banach space. The following holds true for some 
$0<\alpha\leq 1$ depending only on $X$:
\par
Let
$S$ be a subset of $X$ and
$f:S\to\Bbb R$ be a function that is uniformly continuous on bounded sets
of $S$. Then for every $r>0$ and
$\varepsilon>0$ there exists a $\Delta$-convex function 
$F_{r,\varepsilon}:X\to \Bbb R$ 
satisfying the following conditions:
\romanitem{(i)} $\inf_S f=\inf_X F_{r,\varepsilon}$ and 
$\arginf{f}=\arginf{F_{r,\varepsilon}}$,
\romanitem{(ii)} $F_{r,\varepsilon}\in \holder{\alpha}(X)$, 
for some $0<\alpha\leq 1$, and
\romanitem{(iii)} $f(x)-\varepsilon\leq F_{r,\varepsilon}(x)\leq f(x)$ 
for every 
$x\in S\cap B_X(r)$.
\endproclaim
\par
\demo{Proof of the Corollary 4}
By the same argument as above, let $\norm$ be an equivalent norm on $X$ with
modulus of smoothness $1+\alpha$ (for some $0<\alpha\leq 1$).
Consider the following simple extension of $f$: 
$$
F(x):=
\left\{ f(x) \hfill\hbox{for $x\in S$}
\atop
+\infty \quad \hbox{otherwise.}\right.
$$ 
Notice that $\arginf{F}=\arginf{f}\subset S$.
It is not hard to see using {\bf Proposition 8}(i) and 
the proof of {\bf Proposition 6}(v) 
that
the sequence $(\Delta_n^{1+\alpha}{F})_{n\in\Bbb N}$, which satisfies 
(i) and (ii) of {\bf Theorem 1},
also converges uniformly on bounded sets of $S$ to $f$.
\cqd
\enddemo
\par
The proof of {\bf Theorem 1} will be done in a general scheme involving 
two main 
steps. First, we explain an {\it extended inf-convolution} formula 
that gives us a 
standard way to approximate functions on $X$. Then, we develop
some convexity techniques
in order to get smooth $\Delta$-convex functions between the 
functions given by the extended inf-convolution formula.
\par
%
%%%%%%%%
\head 2. 
The extended inf-convolution
\endhead
%%%%%%%%  
\par
In this section we explain the convergence results we need in the proof of 
{\bf Theorem 1}. First, 
we introduce the definition of {\it extended inf-convolution}. This definition
generalizes the classical one of {\it inf-convolution} (see \cite{St$_1$}
for a general survey of the subject) and will be an important tool 
in our work.
\par
\definition{Definition}
For any application $K:X\times X\to \extR$ and any function $f:X\to \extR$ 
we define the {\sl extended inf-convolution} of $f$ by $K$ as the function
$$
\infconv{f}{K}(x):= \inf_{y\in X} \Big\{ f(y) +  K(x,y) \Big\},
\quad
x\in X.
$$
$K$ will be called the {\sl kernel} of the extended inf-convolution.
\enddefinition
\par
\example{Example}
If for $g:X\to \extR$ we consider the kernel $K_g(x,y):=g(x-y)$, then the
extended inf-convolution $\infconv{f}{K_g}$ is nothing else but the classical
inf-convolution $\infconv{f}{g}$. 
\endexample
\par
Before the statement of the main result of this section, 
we need to define some natural properties of kernels.
\par
\definition{Definition}
A kernel $K$ is {\sl pointwise separating} 
if for every $x_0\in X$
and every $\delta>0$ there exists $C_{x_0,\delta}>0$ such that 
$K(x_0,y)\geq C_{x_0,\delta}$ whenever $\|x_0 -y \|\geq \delta$.
\par
A kernel $K$ is called {\sl uniformly separating on bounded sets} 
if for all $r>0$ and $\delta>0$ there exists $C_{r,\delta}>0$ so that 
$K(x,y)\geq C_{r,\delta}$ provided $\|x\|\leq r$ and $\|x-y\|\geq\delta$.
\par
A kernel $K$ is {\sl uniformly separating} if for every 
$\delta>0$  there is some $\beta_{\delta}>0$
in such a way that 
$K(x,y)\geq\beta_{\delta}\|x-y\|$
whenever $\|x-y\|\geq \delta$.
\enddefinition
\par
\definition{Definition}
Given a function $f:X\to \extR$ and a kernel $K$, we define the following 
sequences of functions:
$$
I_{K,n} f := \infconv{f}{nK}
\quad \text {and} \quad
S_{K,n} f := - \infconv{-f}{nK}
\quad 
(n\in\Bbb N).
$$
\enddefinition
\par
\remark{Remark} 
For any Hilbert norm $\norm$  
consider the kernel 
$K_L(x,y)=\|x-y\|^2$.
Then, with our notation 
the sequence
$
\Big( S_{K_L,m} \big(I_{K_L,n} f\big)\Big)_{m>n}
$ 
denotes
the {\sl Lasry-Lions approximates} of $f$ related to the norm $\norm$.
\par
\endremark
\remark{Remark}
Note that the Lasry-Lions approximates commutes with translations in the same
way as the classical inf-convolution also does. This is a consequence of the
following property of the kernel: 
$K_L(x-a,y)=K_L(x,y+a)$ (for all $x$, $y$ and $a$). 
However, the problem of 
regularizing (not necessarily convex) functions in a non-Hilbert space
leads naturally to more general kernels which do not yield 
translation-invariant
approximates.
\endremark
\par
The next facts are easy to check.
\par
\proclaim {Facts 5} 
Let $f:X\to\extR$ be a function.
%
%\roster
%
\romanitem{{\bf 1}} For $x\in X$,
$$
\displaylines{
I_{K,n}f 
=
\inf_{y\in X} \Big\{f(y)+nK(x,y)\Big\},
\cr
S_{K,n} f (x) 
= 
- I_{K,n} (-f) (x)
= 
\sup_{y\in X} \Big\{ f(y) - nK(x,y) \Big\}.
}
$$
\romanitem{{\bf 2}} Let $C$ be a constant. Then $I_{K,n} (f+C)=I_{K,n}f +C$,
for any $n$.
\romanitem{{\bf 3}} Suppose that the kernel $K$ is positive
({\it i.e.}, $K(x,y)\geq 0$ for all $x,y \in X$) then
%
%\roster
\romanitem{(i)} $\big(I_{K,n} f\big)_{n\in\Bbb N}$ is an increasing sequence 
of functions bounded below by $\inf f$.
\romanitem{(ii)} If $f\leq g$, then $I_{K,n} f \leq I_{K,n} g$ for any $n$.
\romanitem{(iii)}$I_{K,m}\big(I_{K,n}f\big) \leq I_{K,m}\big(I_{K,m}f\big)$,
for any $m>n$.
%\endroster
%\endroster
\endproclaim
\par
We now proceed to state and prove a technical proposition which is 
the main result of this section.
\par
\proclaim{Proposition 6}
Let $K:X\times X\to \Bbb R$ a kernel satisfying the following conditions:
\par
%
%\roster
%
\romanitem{(1)} $K$ is {\sl positive}  and $K(x,x)=0$ for all $x\in X$,
\romanitem{(2)} $K$ is {\sl symmetric} ({\it i.e.}, $K(x,y)=K(y,x)$ for all 
$x,y\in X$),
\romanitem{(3)} $K(x,y) @>>y\to\infty> +\infty$ uniformly on bounded sets,
\romanitem{(4)} $K$ is uniformly continuous (resp. Lipschitz continuous) 
on bounded sets and
\romanitem{(5)} $K$ is pointwise separating.
%
%\endroster
%
\par
Then for every proper \lsc bounded below function $f:X\to \extR$ 
the following statements hold:
\par
%
%\roster
%
\romanitem{(i)} $I_{K,n} f  \leq S_{K,n} \big( I_{K,n} f\big) \leq f$.
\romanitem{(ii)} $\inf I_{K,n} f = \inf f$ and 
$\arginf{I_{K,n} f} = \arginf{f}$.
\romanitem{(iii)} $I_{K,n} f$ is uniformly continuous 
(resp. Lipschitz continuous) on bounded sets.
\romanitem{(iv)} $I_{K,n}\big(I_{K,n} f\big) @>\tau_p>n\to\infty>f$ and 
$I_{K,n}\big(I_{K,n} f\big) @>\tau_K>n\to\infty>f$ when $f$ is continuous.
\par
If in addition $K$ is uniformly separating on bounded sets then
\romanitem{(v)} $I_{K,n}\big(I_{K,n} f\big) @>\tau_b>n\to\infty>f$ 
whenever $f$ is uniformly continuous on bounded sets.
\par
Finally, when $K$ is uniformly separating one has
\romanitem{(vi)} $I_{K,n}\big(I_{K,n} f\big) @>\tau_u>n\to\infty>f$ 
provided $f$ is uniformly continuous on $X$ (not necessarily bounded below).
%\endroster
%
\endproclaim
\par
\remark{Remark}
The sequence of functions $I_{K,n}\big(I_{K,n} f\big)$  plays an important
auxiliary r\^ole in this work; namely, it provides a lower bound for
the sequence $(\Delta_{K,n}f)_{n\in\Bbb N}$ in {\bf Proposition 8}(i). 
\endremark
\par
\demo{Proof of the Proposition 6}
\par\noindent
{\it (i)} Since $K(x,x)=0$ we get that $I_{K,n} f \leq f$ (take $y=x$ in the 
infimal definition of $I_{K,n}f$ at any point $x\in X$). 
Therefore we deduce that
$$
S_{K,n}(I_{K,n}) f =- I_{K,n} (-I_{K,n} f)\geq I_{K,n} f.
$$
To see the other inequality, notice that from {\bf Fact 5-1} 
we obtain for $x\in X$ the expression
$$
S_{K,n} (I_{K,n} f) (x) = \sup_{y\in X} \inf_{z\in X} 
\Big\{ f(z) + n\big(K(y,z) - K(x,y)\big) \Big\}.
\eqno(\nfn)
$$ 
For some  fixed $x$, if we take $z=x$  in \pfn{0} 
we conclude from the symmetry of $K$ that $S_{K,n}(I_{K,n}f) (x) \leq f(x)$.  
\par\noindent
(ii) From {\bf (i)} and {\bf Fact 5-1}(i) we have      
$\inf I_{K,n} f = \inf f$ and $\arginf{f} \subseteq \arginf{I_{K,n} f}$.
Consider any minimum $x_0\in X$ of $I_{K,n}f$. Then, there exists a sequence
$(y_k)_{k\in\Bbb N}\subset X$ so that 
$$
\inf f = I_{K,n}f (x_0)\leq f(y_k) + nK(x_0,y_k) @>>k\to\infty>\inf f.
\eqno(\nfn)
$$
Hence, since $K$ is positive it follows from \pfn{0} that
$$
\lim_{k\to\infty} f(y_k) = \inf f 
\text { and }
\lim_{k\to\infty} K(x_0,y_k)=0.
\eqno(\nfn)
$$ 
But $K$ is pointwise separating, so the second part of \pfn{0} 
implies that
$y_k @>>>x_0$. Using the lower-semicontinuity of $f$ and the first part of 
\pfn{0} we conclude that
$$
\inf f\leq f(x_0)\leq \lim_{k\to\infty} f(y_k) = \inf f.
$$
and this proves assertion (ii).
\par
Before proceeding with the rest of the proof, we set up the following useful
definition:
$$
\Omega_n (x):=\big\{ y\in X : f(y)+nK(x,y) \leq I_{K,n}f (x)+1\big\} \quad
(x\in X, n\in \Bbb N) 
\eqno(\nfn)
$$
With these notations, we remark that for $n\in N$ and $x\in X$
$$
I_{K,n}f (x)=\inf_{y\in\Omega_n(x)} \big\{ f(y)+nK(x,y)\big\}
\geq
\inf_{\Omega_n(x)} f
\eqno(\nfn)
$$
(the last inequality coming from the positivity of $K$).
\par
It is clear from \pfn{0} that the behaviour of $I_{K,n}f$ 
is directly linked with the size of the sets 
$\big\{\Omega_n(x)\big\}_{x\in X}$.
We shall see that the growth condition (3) ensures that the sets  
$\Omega_n(x)$ are not arbitrarily big when $x$ 
runs on bounded sets of $X$. 
More precisely, we claim the following.
\par
\proclaim{Claim 6.1}
For any $r>0$, the set 
$\Omega_r:=\bigcup_{n\in\Bbb N} \bigcup_{\|x\|\leq r} \Omega_n(x)$ 
is bounded.
\endproclaim
\par
The proof of this claim is based on the next simple fact.
\par
\proclaim{Fact 6.2}
For any $r>0$,
$
\sup\Big\{ \frac {I_{K,n}f (x)}{n}:x\in B_X(r),\, n\in\Bbb N\Big\}
:=M_r <+\infty
$.
\endproclaim
\par
\demo{Proof of the Fact 6.2}
Since $f$ is proper, take $y_0$ such that $f(y_0)\leq \inf f +1< +\infty$.
Then by definition of $I_{K,n}f$ it follows that for any $x\in X$ 
$$
\frac {I_{K,n}f(x)}{n} 
\leq 
\frac{f(y_0)}{n}+K(x,y_0)
\leq 
\inf f + 1 +\sup \big\{ K(x,y_0): x\in B_{X}(r)\big\},
$$
and this expression is bounded above on bounded sets because $K$ 
is uniformly continuous (or Lipschitz continuous) on bounded sets.
The proof of {\bf Fact 6.2} is finished.
\enddemo
\par
\demo{Proof of the Claim 6.1}
For $r_0>0$, let $M_{r_0}>0$ be the upper bound defined in {\bf Fact 6.2}. 
Thus, for any $x\in B_X(r_0)$ and $n\in\Bbb N$ if $y\in\Omega_n(x)$ 
it follows from the definition of $\Omega_n(x)$, given in 
\pfn{-1}, that
$$
K(x,y)\leq \frac {1}{n} \Big( I_{K,n}f(x)+1-f(y)\Big)\leq M_{r_0}+1-\inf f.
\eqno(\nfn)
$$
But the growth condition on $K$ given by (3) implies that the set of $y$ 
satisfying \pfn{0} is uniformly bounded for $x\in B_X(r_0)$. The proof of
{\bf Claim 6.1} is done.
\enddemo
\par
We can now continue with the proof of Proposition 6.
\par
\noindent
(iii) Suppose the kernel $K$ is Lipschitz continuous on bounded sets
(the proof for the uniformly continuous case is practically the same). 
For $r_0>0$ take $x,x'\in B_X(r_0)$ and let $L_{K,r_0}$ be the Lipschitz 
constant of $K$ on $B_X(r_0)\times\Omega_{r_0}$ ($\Omega_{r_0}$ being 
bounded by {\bf Claim 6.1}). 
Using the equality of \pfn{-1}
we can construct a sequence $(y_k)_{k\in\Bbb N}\subset\Omega_{r_0}$ 
in such a way that for every $k\in\Bbb N$ one has
$
f(y_{k})+nK(x',y_{k})\leq I_{K,n} f(x') + \frac {1}{k}
%\eqno(\nfn)
$.
Therefore, we obtain
$$
\eqalign{
I_{K,n}f (x') - I_{K,n} f(x) 
& \leq 
f(y_k) +nK(x',y_k) - f(y_k) - nK(x,y_k)+\frac{1}{k}
\cr
&\leq 
nL_{K,r_0}\|x'-x\| + \frac{1}{k}
@>>k\to\infty>nL_{K,r_0}\|x'-x\|.
}
$$
This concludes the proof of (iii).
\par
We first prove (iv), (v) and (vi) for 
$(I_{K,n} f)_n$ instead of
$\big(I_{K,n}(I_{K,n} f)\big)_n$. We will complete the proof afterwards.
\par
\noindent
(iv') Fix $x_0\in X$. If 
$\lim_{n\to\infty} I_{K,n} f(x_0) = \sup_{n} I_{K,n} f(x_0) = +\infty$
then by {\bf (i)} one has $f(x_0)=+\infty$ and the result holds. 
Thus, suppose that
$I_{x_0}:= \lim_{n} I_{K,n} f(x_0)<+\infty$. By the infimal definition of 
$I_{K,n}f$ at $x_0$, we can choose a sequence 
$(y_n)_{n\in\Bbb N} \subset X$ such that
$$
I_{K,n} f (x_0)
\leq 
f(y_n) + nK(x_0,y_n) 
\leq
I_{K,n} f(x_0) + \frac {1}{n}
@>>n\to\infty> I_{x_0}
\eqno(\nfn)
$$
Hence, from \pfn{0} it follows for $n\in\Bbb N$ that
$$
K(x_0,y_n)
\leq 
\frac {1}{n}\big( I_{K,n} f (x_0) - f(y_n) \big) +\frac {1}{n^2}
\leq 
\frac {1}{n}\big( I_{x_0} - \inf f \big) +\frac{1}{n^2}
@>>n\to\infty> 0.
\eqno(\nfn)
$$
But $K$ is pointwise separating, so we have from \pfn{0}
that $(y_n)_n$ is norm converging to $x_0$. 
Using the lower-semicontinuity of $f$, 
the positivity of $K$ in \pfn{-1} and {\bf (i)}, we get that
$$
f(x_0)
\leq
\liminf_{n\to\infty} f(y_n)\leq I_{x_0}\leq f(x_0).
$$
\par
If $f$ is continuous, since by {\bf Fact 5-3}(i)  and {\bf (iii)}
$(I_{K,n}f)_n$ is an increasing sequence of continuous functions, 
Dini's Theorem tell us that the pointwise convergence of 
$(I_{K,n}f)_n$ to $f$ is actually uniform on compact sets.
\par
\noindent
(v') Let $f$ be an uniformly continuous function on bounded sets and 
$O_{r_0}$ be the oscillation of $f$ on the set 
$B_X(r_0)\cup\Omega_{r_0}$, for some fixed $r_0>0$. 
Then, for any $n\in\Bbb N$,
$x\in B_X(r_0)$ and $y\in\Omega_n(x)$ after the first inequality 
of \pfn{-2} and 
{\bf (i)} we have that
$$
K(x,y)
\leq 
\frac {1}{n} \big( I_{K,n} f(x)+1-f(y)\big)
\leq
\frac {1}{n} \big( f(x)-f(y)+1\big)
\leq
\frac {1}{n} (O_{r_0} +1)
@>>n\to\infty>
0.
\eqno(\nfn)
$$
\par
Suppose that $K$ is uniformly separating on bounded sets .
Then, a direct consequence 
of \pfn{0} is that $\lim_n\diam(\Omega_n(x))=0$ uniformly on $B_X(r_0)$. 
Therefore, it follows from {\bf (i)}, \pfn{-4} and the uniform 
continuity of $f$ on $B_X(r_0)$
that
$$
f(x)
\geq 
\lim_{n\to\infty} I_{K,n}f(x)
\geq 
\lim_{n\to\infty} \inf_{\Omega_n(x)}f
@>>n\to\infty> 
f(x)
\eqno(\nfn)
$$
uniformly on $x\in B_X(r_0)$.
\par
\noindent
(vi') Suppose that $f$ is uniformly continuous on $X$. 
Then $f$ satisfy the following fact 
(whose simple proof is left as an exercise to the reader):
$$
\text{there exists } \alpha>0 \text { such that }
f(x)-f(y) \leq \max \{ 1, \alpha\|x-y\| \}
\text{ for all } x,y \in X.
\eqno(\nfn)
$$
\par
Then, in 
the same way as in \pfn{-2} before, using this time \pfn{0}, we deduce 
that for $n\in\Bbb N$, $x\in X$ and any $y\in\Omega_n(x)$ 
$$
K(x,y)
%\leq
%\frac {1}{n} \big(I_{K,n}f(x)+1-f(y)\big)
\leq
\frac {1}{n} (f(x)-f(y)+1)
\leq
\max\Big\{\frac {1}{n},\frac {\alpha}{n}\|x-y\| \Big\}+\frac{1}{n}.
\eqno(\nfn)
$$
\par
For $1>\delta>0$, since $K$ is uniformly separating there is some 
$\beta_\delta>0$ so that from \pfn{0} we deduce for 
$x\in X$ and $y\in\Omega_n(x)$ that
$$
\|x-y\|
\leq
\max 
\Big\{ 
\frac{1}{n\beta_\delta}, \frac{\alpha}{n\beta_\delta}\|x-y\|
\Big\}
+\frac{1}{n\beta_\delta}
\text{ whenever } \|x-y\|>\delta.
\eqno(\nfn)
$$
Hence, taking $n$ big so that 
$
\max\big\{\frac{2}{n\beta_\delta},\frac{2\alpha}{n\beta_\delta}\big\}
\leq\delta<1
$,
\pfn{0} shows for every $x\in X$ that 
$\diam \big(\Omega_n(x)\big)\leq 2\delta$.
\par
That is, we have shown that $\diam\big(\Omega_n(x)\big)\to 0$ uniformly on 
$x\in X$.
Therefore, as $f$ is uniformly continuous on $X$ we can 
repeat the same reasonings of \pfn{-3} to conclude that $(I_{K,n} f)_n$ 
converges to $f$ uniformly on $X$.
\par
(iv) and (v) are straightforward corollaries of (iv') and (v') 
if we remark the following.
\par
Suppose that for $\varepsilon>0$ there exists $n_0\in\Bbb N$ so that 
$%$
f -\frac{\varepsilon}{2}
\leq
I_{K,n_0} f
%\leq 
%f
%\eqno(\nfn)
$%$
\ on some set $S$ 
($S$ being a singleton, or a compact set or a bounded set of $X$). 
By {\bf Fact 5-3}(i) 
and {\bf (iii)}, 
we can then apply (iv') (or  (v')) to the bounded below, 
uniformly continuous function $I_{K,n_0}f$ to obtain $m>n_0$ such that
$%$
I_{K,n_0}f - \frac{\varepsilon}{2}
\leq
I_{K,m}\big(I_{K,n_0} f\big)
%\leq 
%I_{K,n_0}f
%\eqno(\nfn)
$%$
\ on the same $S$. Thus, by {\bf Fact 5-3}(iii) and {\bf (i)} 
it follows that
$$
f-\varepsilon 
\leq 
I_{K,n_0}f 
-
\frac {\varepsilon}{2}
\leq
I_{K,m}\big(I_{K,n_0} f\big)
\leq
I_{K,m}\big(I_{K,m} f\big)
\leq
f
\text { on } S.
$$
\par
(vi) is also easily deduced from (vi') through the following argument. 
\par
If $f-{\varepsilon}\leq I_{K,n}f\leq f$,
for some $\varepsilon>0$ and $n\in\Bbb N$,
then applying {\bf Facts 5-2} and {\bf 5-3}(ii) we get that
$$
f-2\varepsilon 
\leq 
I_{K,n}f - {\varepsilon}
=
I_{K,n}(f - {\varepsilon})
\leq
I_{K,n}\big(I_{K,n}f\big)
\leq
I_{K,n} f
\leq
f.
\cqd
$$
\enddemo
\par
\remark{Remark}
With the above techniques it is not difficult to check
that $\big(I_{K,n}(I_{K,n} f)\big)_n$ converges to $f$ 
for the {\sl epigraphical distance}
(see \cite{AW} for the definition). We refer to the proof of 
{Lemma 3}(v) in \cite{St$_2$} for details.
\endremark
\par
%
%%%%%%%%
\head 3. 
Convexity techniques and smoothness results
\endhead
%%%%%%%%
\par
In this section we shall show a procedure to obtain smooth functions 
from the operators $I_{K,n}(\cdot)$ and $S_{K,n} (\cdot)$. 
We will need to impose some additional conditions of
convexity and smoothness  
on the kernel $K$ to achieve the smooth regularization. 
The interesting feature of these convexity arguments is the 
preservation of the approximating properties obtained in
the previous section.
\par
The main tool we shall use to get smooth regularization is 
explained in the next theorem. It deals with the smooth 
properties inherited by the convex envelop of a ``somehow'' 
smooth function.
\par
\proclaim{Theorem 7}
Let $c:X\to {\Bbb R}$ be a differentiable function, 
and $d:X\to {\Bbb R}$ be a convex function. Denote by $h$
their difference $h:=c-d$ and assume that $\coenv{h}$ makes sense.  
Then the following statements are fulfiled.
\par
%
%\roster
%
\romanitem{(i)} If $c\in\unifder(X)$ (resp. $c\in\holder{\alpha}(X)$, for 
some $0<\alpha\leq 1$) then
$\coenv h\in\unifder(X)$ (resp. $\coenv h\in\holder{\alpha}(X)$). 
\romanitem{(ii)} If $c\in\bdunifder(X)$ (resp. $c\in\bdholder{\alpha}(X)$, 
for some $0<\alpha\leq 1$) and $h$ 
is uniformly continuous
on bounded sets and {\sl strongly coercive} ({\it i.e.},
$
\lim\limits_{x\to\infty}\frac {h(x)}{\|x\|}=+\infty
$)
then $\coenv h\in\bdunifder(X)$ (resp. $\coenv h\in\bdholder{\alpha}(X)$). 
%
%\endroster
%
\endproclaim
\par
\remark{Remark} A proof for the finite dimensional version of
{\bf Theorem 7}(ii) with $d\equiv 0$ can be found in \cite{GR}. 
Our more general proof does not require local compactness and 
relies upon ideas of the work \cite{Fa}. 
The fact that the convex envelope of a smooth function $c$ ``perturbed'' by
a non-smooth concave function $-d$ is still smooth will be crucial later
(namely,
when we check the smoothness of the sequence $(\Delta_{K,n}f)_{n\in\Bbb N}$
in {\bf Proposition 8}).
\par
Notice that the uniform continuity hypothesis on the derivative of $c$ 
cannot be weakened
in the infinite dimensional setting. There are 
bounded below ${\Cal{C}}^{\infty}$-differentiable functions
on $\ell_2$ whose convex envelope is not even G\^ateaux differentiable 
(see Example II.5.6(a) in \cite{DGZ}).
\endremark
\par
\demo{Proof of the Theorem 7}
Denote by $\nu:=\coenv{h}=\coenv(c-d)$.
\par
\noindent
(i) Suppose that $c\in\holder{\alpha}(X)$ (the proof for the other case is 
similar). Since $\nu$ is convex, a necessary and sufficient condition for
$\nu\in\holder{\alpha}(X)$ is that for every $x,\,y\in X$ one has
$$
\nu(x+y)+\nu(x-y)-2\nu(x)
\leq
L\|y\|^{1+\alpha}, 
\text{ for some } L>0.
\eqno(\nfn)
$$
(see Lemma V.3.5 of \cite{DGZ}). We shall check this 
condition for $\nu$.
\par
For $\varepsilon>0$ and $x\in X$, by the expression of the convex envelope of
a function given in \pfn{-14}, we can choose $x_1,\dots,x_n\in X$ and
$\lambda_1,\dots,\lambda_n>0$ 
so that
$$
\sum_{i=1}^n \lambda_i =1,\
\sum_{i=1}^n \lambda_i x_i=x
\text { and }
\sum_{i=1}^n \lambda_i h(x_i)\leq \nu(x) + \frac{\varepsilon}{2}.
\eqno(\nfn)
$$
Note that from the two first parts of \pfn{0} we also have 
$$
x\pm y
=
\bigg(\sum_{i=1}^n \lambda_i x_i\bigg) 
\pm 
\bigg(\sum_{i=1}^n \lambda_iy\bigg)
=
\sum_{i=1}^n \lambda_i(x_i\pm y).
\eqno(\nfn)
$$
Thus, it follows from \pfn{-16} and \pfn{0} that 
$$
\nu(x\pm y) \leq \sum_{i=1}^n \lambda_i h(x_i\pm y).
\eqno(\nfn)
$$
Let $L>0$ be the $\alpha$-H\"older continuity constant 
of the derivative of $c$.
Putting together the last part of \pfn{-2}, \pfn{0} 
and using the convexity of $d$, we get
$$
\displaylines{
\nu(x+y)+\nu(x-y)-2\nu(x)
\leq
\sum_{i=1}^n \lambda_i h(x_i+y) + \sum_{i=1}^n \lambda_i h(x_i-y)-
2\sum_{i=1}^n \lambda_i h(x_i) +\varepsilon=
\cr
\sum_{i=1}^n\!\lambda_i\Big(c(x_i+y)+c(x_i-y)-2c(x_i)\Big)
+
%\sum_{i=1}^n\!\lambda_i\Big(d(x_i+y)+d(x_i-y)-2d(x)\Big)+\varepsilon
2\sum_{i=1}^n\!\lambda_i\bigg(d(x_i)-\frac{d(x_i+y)+d(x_i-y)}{2}\bigg)
+\varepsilon
\leq
\cr
\sum_{i=1}^n\!\lambda_i\Big(c(x_i+y)-c(x_i)+c(x_i-y)-c(x_i)\Big)+\varepsilon 
\leq\!%
\sum_{i=1}^n\!\lambda_i 2^{\alpha}L\|y\|^{1+\alpha} +\varepsilon 
\!=\!
2^{\alpha}L\|y\|^{1+\alpha} +\varepsilon,
}
$$
because $c\in\holder{\alpha}(X)$ (and therefore satisfy \pfn{-3} for
$L'=2^{\alpha}L$). As $\varepsilon$ is arbitrary, the condition \pfn{-3}
holds for $\nu$.
\par
\noindent
(ii) can be proved reproducing the same lines as before, bearing in mind
that the lack of uniformity for the derivative of $c$ on $X$ can be replaced
by the next ``localization'' property of the convex envelope of a strongly
coercive function $h$.
\par
\proclaim{Claim 7.1}
Let $h:X\to\Bbb R$ be a function which is 
uniformly continuous on bounded sets 
and strongly coercive. Then for every $r>0$ there
exists $\rho_r>0$ so that for all $\|x\|\leq r$ one has
$$
\coenv{h}(x)=
\inf
\bigg\{ \sum_{i=1}^n \lambda_{i} h(x_i):
(x_i)_{i=1}^n\subset B_X(\rho_r),\,
\lambda_i>0,\,
\sum_{i=1}^n \lambda_i=1,\,
x=\sum_{i=1}^n \lambda_i x_i
\bigg\}.
$$
\endproclaim
\demo{Proof of the Claim 7.1}
First, note that under the hypothesis of {\bf Claim 7.1}, $h$ is bounded
below, so that $\coenv{h}$ makes sense.
Fix $r_0>0$ and let $m_{r_0}$  be the infimum of $h$ on $X$ and $M_{r_0}$
be the supremum of $h$ on $B_X(r_0+1)$ 
($M_{r_0}<+\infty$, because of the uniform continuity
of $f$ on $B_X(r_0+1)$).
Consider the following family of
hyperplanes:
$$
{\Cal H}_{r_0}
:=
\bigcup_{\textstyle {x\in B_X(r_0)\atop v\in B_X,\, v^*\in B_{X^*}}}
\Big\{ 
H_{x,v}(z)= m_{r_0} +\big(h(x+v)-m_{r_0}\big)v^*(z-x) :
v^*(v)=1 
\Big\}.
\eqno(\nfn)
$$
\par
Notice that for $\|x\|\leq r_0$ and $v\in B_X$ we have the following
$$
H_{x,v} (x)=m_{r_0}
\leq
\coenv{h}(x)
\text{ and }
H_{x,v} (x+v)=h(x+v)\leq M_{r_0}.
\eqno(\nfn)
$$
\par
Since $h$ is strongly coercive, we get from \pfn{-1} that
$$
\sup_{H\in{{\Cal H}_{r_0}}} H(z)
\leq
m_{r_0} + (M_{r_0}-m_{r_0})(\|z\|+r_0)
<
h(z)
\text{ provided }
\|z\|>\rho_{r_0},
\eqno(\nfn)
$$
for some $\rho_{r_0}>0$. Let us show that $\rho_{r_0}$ satisfy 
the conclusion of the claim.
\par
The strategy is to replace any convex combination that appears in the 
definition
of the convex envelope \pfn{-20} by another smaller convex combination
with ``uniformly bounded vertices''. This idea is formally explained in
the next fact. 
\par
\proclaim{Fact 7.2}
For  $\|x\|\leq r_0$, 
consider any finite convex combination 
($x_1,\dots,x_n\in X$, $\lambda_1,\dots,\lambda_n>0$ and $\sum_{i=1}^n
\lambda_i=1$) such that $\sum_{i=1}^n \lambda_i x_i=x$. 
If $\|x_{i_0}\|>\rho_{r_0}$, 
for some $1\leq i_0\leq n$,
then there exists $x_{i_0}'\in B_X(\rho_{r_0})$
and $\lambda_1',\dots,\lambda_n'>0$ so that 
$
\sum_{i=1}^n \lambda_i'=1 
$,
$$ 
\sum_{\textstyle {i=1\atop i\not=i_0}}^n 
\lambda_i'x_i + \lambda_{i_0}' x_{i_0}'=x
\text { and }
\sum_{\textstyle {i=1\atop i\not=i_0}}^n 
\lambda_i'h(x_i) 
+ \lambda_{i_0}'h(x_{i_0}')
\leq
\sum_{i=1}^n \lambda_i h(x_i).
$$
\endproclaim
\demo{Proof of the Fact 7.2}
For simplicity, take $i_0=1$. 
Since $\|x_1\|>\rho_{r_0}$, it follows from
\pfn{0} that $H_{x,v_{x_1}}(x_1)< h(x_1)$,
where we take $v_{x_1}:=\frac{x_1-x}{\|x_1-x\|}$. But by the first part of 
\pfn{-1} we also have
$H_{x,v_{x_1}}(x)=m_{r_0}\leq\sum_{i=1}^n \lambda_i h(x_i)$. Hence, the
segment 
$I_{x,x_1}
:=
\Big[ \big(x,\sum_{i=1}^n \lambda_i h(x_i)\big), \big(x_1,h(x_1)\big)\Big]
\subset X\times \Bbb R
$
belongs to the upper half-space define by $H_{x,v_{x_1}}$. 
Therefore, the equality in the second part of \pfn{-1} implies that
$$ 
I_{x,x_1}
\cap 
\big\{ (x+v_{x_1}, t): t\in\Bbb R\big\} 
=
\big\{(x+v_{x_1}, s)\big\}
\text{ where } 
h(x+v_{x_1})\leq s.
\eqno(\nfn)
$$
\par
If we define 
$x':=\sum_{i>1}^n \frac{\lambda_i}{1-\lambda_1}x_i$ (so that
we have $x=\lambda_1 x_1 +(1-\lambda_1)x'$), using barycentric
coordinates on the segment 
$
\Big[
\big(x',\sum_{i>1}^n \frac{\lambda_i}{1-\lambda_1}h(x_i)\big),
\big(x_1,h(x_1)\big)\Big]
$
we can compute some $\mu\geq 0$ 
in such a way that
$$
\eqalign{
\Bigg(x,\sum_{i=1}^n \lambda_i h(x_i)\Bigg)
&=
\mu\Bigg(x',\sum_{i>1}^n \frac{\lambda_i}{1-\lambda_1}h(x_i)\Bigg) 
+ 
(1-\mu) (x+v_{x_1}, s)
=
\cr
&\Bigg(
\sum_{i>1}^n\frac{\mu\lambda_i}{1-\lambda_1}x_i,
\sum_{i>1}^n\frac{\mu\lambda_i}{1-\lambda_1}h(x_i)
\Bigg)
+ 
(1-\mu)(x+v_{x_1},s).
}
\eqno(\nfn)
$$
\par
But \pfn{-1} and \pfn{0} together give that
$
x
=
\sum_{i>1}^n\frac{\mu\lambda_i}{1-\lambda_1}x_i
+
(1-\mu)(x+v_{x_1})
$
and 
$$
\sum_{i>1}^n\frac{\mu\lambda_i}{1-\lambda_1}h(x_i)
+
(1-\mu) h(x+v_{x_1})
\leq
\sum_{i>1}^n\frac{\mu\lambda_i}{1-\lambda_1}h(x_i)
+
(1-\mu)s
=
\sum_{i=1}^n \lambda_i h(x_i).
$$
\par
This concludes the proof of {\bf Fact 7.2}. 
\enddemo
Then the proof of {\bf Claim 7.1} is done and, therefore, 
{\bf Theorem 7} is proved.
\cqd
\enddemo
\enddemo
\par
\remark{Remark}
{\bf Claim 7.1} is false for functions $h$ failing the strong coerciveness
condition
$\liminf\limits_{x\to\infty} \frac {h(x)}{\|x\|}=+\infty$. For instance,
consider $h:\Bbb R\to\Bbb R$ defined by $h(x)=\sqrt{|x|}$.
\endremark
\par
{\bf Theorem 7} can be applied as an useful tool to regularize functions
on infinite dimensional Banach spaces. Our next proposition, 
which provides the
smoothness assertions we need for proving {\bf Theorem 1}, is a good
example of this feature. We keep the notation used in {\tensmc Section 2}.
\par
\proclaim{Proposition 8}
Let $K:X\times X\to \Bbb R$ be a kernel satisfying the following conditions:
\par
\romanitem{(1)} $K$ is {\sl positive}  and $K(x,x)=0$ for all $x\in X$,
\romanitem{(2)} $K$ is {\sl symmetric}, 
\romanitem{(3)} $K(x,y) @>>y\to\infty> +\infty$ uniformly on bounded sets,
\romanitem{(4)} $K$ is uniformly continuous
on bounded sets and
\romanitem{(5)} $K(x,y)= c_K(x)-d_K(x,y)$ 
where $d_K(\cdot,y)$ is a lower semi-continuous convex function, 
for all $y\in X$.
\par
Let $f:X\to\extR$ a proper lower semi-continuous function 
and consider the sequence of 
$\Delta$-convex
functions $\Delta_{K,n}f$ defined as follows
$$
\Delta_{K,n}f:= \coenv{(I_{K,n}f +nc_K)} - nc_K
\quad (n\in\Bbb N).
$$
Then the following assertions hold.
\par
\romanitem{(i)} $I_{K,n}\big(I_{K,n}f\big) \leq \Delta_{K,n}f\leq f$.
\romanitem{(ii)} If $c_K\in\unifder(X)$ 
(resp. $c_K\in\holder{\alpha}(X)$, for some $0<\alpha\leq 1$), then 
one has
$(\Delta_{K,n}f)_{n\in\Bbb N}\subset\unifder(X)$
(resp. $(\Delta_{K,n}f)_{n\in\Bbb N}\subset\holder{\alpha}(X)$).
\romanitem{(iii)} If $c_K\in\bdunifder(X)$ 
(resp. $c_K\in\bdholder{\alpha}(X)$, for some $0<\alpha\leq 1$)
and $c_K$ is strongly coercive
then
$(\Delta_{K,n}f)_{n\in\Bbb N}\subset\bdunifder(X)$
(resp. $(\Delta_{K,n}f)_{n\in\Bbb N}\subset\bdholder{\alpha}(X)$)
provided $f$ is bounded below. 
\par
\endproclaim
\par
\remark{Remark}
For every pair of function $f,g$ on $X$, denote by 
$\Deltconv{f}{g}:=\coenv{(f+g)}-g$. Suppose that $\norm$ is 
a Hilbert norm and
consider the kernel $K_L(x-y):=\|x-y\|^2$.
The Lasry-Lions approximates 
of a function $f$ by the norm $\norm$ satisfy the following relation for
$m>n$ (see Proposition 2(i) of \cite{St$_2$})
$$
\Big( S_{K_L,m} \big(I_{K_L,n} f\big)\Big)
=
I_{K_L,m-n}\big(\Deltconv{f}{nc}\big).
$$
Compare with the expression given by {\bf Proposition 8}
$$
\Delta_{K_L,n}f=\Deltconv{\big(I_{{K_L},n}f\big)}{nc}.
$$
\endremark
\par
\demo{Proof of the Proposition 8} 
For any function $g:X\to \extR$, 
denote by 
$$
D_n g(x):=\sup\{g(y)+nd_K(x,y):y\in X\}.$$
Since $d_K(\cdot,y)$
is a \lsc convex function, we have that $D_{n}g$ is \lsc and convex too.
Note that by {\bf Fact 5-1} we also have the following 
decomposition:
$$
S_{K,n}g
=
\sup_{y\in X} \Big\{g(y)-n\big(c_K(x)-d_K(x,y)\big)\Big\}
=
D_n g(x) - nc_K(x).
\eqno(\nfn)
$$
\par
On the other hand, (1) and (2) 
ensure that (i) of {\bf Proposition 6} holds true.
Hence, for all $n\in\Bbb N$ from \pfn{0} we get that
$$
I_{K,n}(I_{K,n} f)
\leq
S_{K,n} \big(I_{K,n}(I_{K,n} f)\big)
=
D_n\big(I_{K,n}I_{K,n} f\big) - nc_K
\leq
I_{K,n} f
\leq  f
\eqno(\nfn)
$$
\par
Now, we make the next simple but crucial observation.
\par
\proclaim{Fact 8.1}
Let $c,d$ and $e$ be three functions such that $d-c\leq e$ and suppose 
that $d$ is \lsc
and convex.
Then we have that $d-c\leq\coenv{(e+c)}-c\leq e$.
\endproclaim
\demo{Proof of the Fact 8.1}
It suffices to 
note that the convexity of $d$ 
implies the equivalent inequality
$
d\leq \coenv{(e+c)}\leq e+c
$. 
\enddemo
\par
Applying {\bf Fact 8.1} to the inequality \pfn{0} we obtain that
$$
I_{K,n}I_{K,n} f
\leq
D_n\big(I_{K,n}I_{K,n} f\big) - nc_K
\leq
\coenv{(I_{K,n}f +nc_K)} - nc_K
=
\Delta_{K,n}f
\leq
I_{K,n} f.
$$
\par
At this point, another important remark turns up. By definition of
$I_{K,n} f$ at 
any point $x\in X$ one has
$$
(I_{K,n} f +nc_K)(x)
=
2nc_K(x)+\inf_{y\in X} \big\{f(y)-d_K(x,y)\big\}
=
2nc_K(x) - D_n (-f)(x),
\eqno(\nfn)
$$
where $D_n(-f)$ is a convex function.
\par
Therefore, if $c_K\in\unifder(X)$ (or $c_K\in\holder{\alpha}(X)$)
{\bf Theorem 7}(i) can be applied to $\coenv{(I_{K,n} f(x) +nc_K(x))}$ 
because of \pfn{0}. This shows {\bf Proposition 8}(ii).
\par
The proof for the case of $c_K\in\bdunifder(X)$ 
(or $c_K\in\bdholder{\alpha}(X)$) can be done
in a similar way from {\bf Theorem 7}(ii). Notice that  
the function $I_{K,n} f +nc_K$ is uniformly continuous on bounded sets
since $K$ satisfies (1)--(4) and therefore {\bf Proposition 6}(iii) holds. 
On the other hand, the strong coerciveness of $c_{K}$ implies
for a bounded below
function $f$ that
$$
\frac {I_{K,n} f(x) +nc_K(x)}{\|x\|}
\geq
\frac {\inf f}{\|x\|} + n\frac {c_K(x)}{\|x\|}
@>>x\to\infty> +\infty.
\cqd
$$
\enddemo
\par
%
%
%
%%%%%%%%
\head 4. 
The proof of the main result
\endhead
%%%%%%%%
%
With the tools shown in {\tensmc Section 2} and {\tensmc Section 3},
we are now ready to prove our main result.
\par
\demo{Proof of the Theorem 1}
Let $f:X\to\extR$ be a bounded below \lsc function.
For $p>1$, we define the following kernel $K_p:X\times X\to \Bbb R$ as
$$
K_p(x,y):=2^{p-1}\|x\|^p + 2^{p-1}\|y\|^p - \|x+y\|^p.
\eqno(\nfn)
$$
Let us first check the following two basic properties of $K_p$.
\par
\noindent
(1) Clearly, $K_p(x,x)=0$ (for all $x\in X$). Also, $K_p$ is positive 
since one has that 
$$
\|x+y\|^p
\leq
\big(\|x\|+\|y\|\big)^p
\leq
2^{p-1}\big(\|x\|^p+\|y\|^p\big).
$$
\par
\noindent
(2) $K_p$ is obviously symmetric.
\par
If we set $c_{K_p}(x):=2^{p-1}\|x\|^p$ and 
$d_{K_p}(x,y):=\|x+y\|^p-2^{p-1}\|y\|^p$ ($x,\,y\in X$),
we have 
$
K_p(\cdot,y)=c_{K_p}(\cdot)-d_{K_p}(\cdot,y).
$
\par
With this notation and the definition of $g_n^p$ given in {\bf Theorem 1},
one has 
$$
(I_{{K_p},n}f+nc_{K_p} )= g_n^p
\quad
\text{for every $n\in\Bbb N$}.
\eqno(\nfn)
$$
Hence, using again the notation of {\bf Theorem 1} 
and {\bf Proposition 8}(i), it follows that
$$
I_{{K_p,}n}\big(I_{{K_p},n}f\big) 
\leq 
\Delta_{K_p,n}f
=
\coenv{g_n^p} - 2^{p-1}n\norm^p
=
\Delta_n^p{f} 
\leq 
f.
\eqno(\nfn)
$$
Therefore, by \pfn{0} 
the statements (i), (iii), (iv) and (v) of {\bf Theorem 1}
hold true if we check that $K_p$ satisfies the assumptions of {
\bf Proposition 6}.
\par
We proceed to show the following growth property of $K_p$, 
that trivially implies the condition (3) of {\bf Proposition 6}.
\par
\proclaim{Claim 1.1}
For any $p>1$ there exists $\gamma_p>0$ and $\eta_p>1$ so that 
$K_p(x,y)\geq \gamma_p \|y\|^p$ whenever $\|y\|\geq\eta_p\|x\|$.
\endproclaim
\demo{Proof of Claim 1.1}
Take $\eta>1$ and $x,y\in X$ such that $\eta\|x\|\leq\|y\|$. After
the computation
$$
\eqalign{
K_p(x,y)
&\geq
\|y\|^p
\bigg(
2^{p-1}\Big|\frac{\|x\|}{\|y\|}\Big|^p +2^{p-1}
-
\Big\|\frac{y}{\|y\|}+\frac{x}{\|y\|}\Big\|^p
\bigg)
\geq\hfill
\cr
&\|y\|^p \Big(2^{p-1}-\Big|1+\frac{\|x\|}{\|y\|}\Big|^p \Big)
\geq
\|y\|^p \Big( 2^{p-1} - \big(1+\frac {1}{\eta}\big)^p\Big).
}
$$
The claim is proved by choosing 
$\eta_p>1$
such that
$\gamma_p:=\big( 2^{p-1} - (1+\frac {1}{\eta_p})^p\big)>0$.
\cqd
\par
\enddemo
\par
It is clear that $K_p$ is Lipschitz continuous on bounded sets. The
next claim takes care of the {\sl separating} properties of $K_p$.
\par
\proclaim{Claim 1.2}
Suppose that the norm $\norm$ is l\.u\.c\. at $x_0\in X$ (resp. UC) then 
for every $\varepsilon>0$ there exists $C_{\varepsilon,x_0}>0$  such that
$K_p(x_0,y)\geq C_{\varepsilon,x_0}\|x_0-y\|^p$ whenever 
$\|x_0-y\|\geq\varepsilon$
(resp. for all $r>0$ and $\varepsilon>0$ there exists $C_{\varepsilon,r}>0$
such that 
$K_p(x_,y)\geq C_{\varepsilon,r}\|x-y\|^p$ provided 
$\|x-y\|\geq\varepsilon$ and $\|x\|\leq r$).
\endproclaim
\par
\demo{Proof of the Claim 1.2}
We only prove the claim under the uniform convexity assumption. 
The proof for the l\.u\.c\.
case is completely similar. We proceed by contradiction.
\par
Suppose that the claim is false. Then by definition of $K_p$ (see \pfn{-2}), 
there are two sequences
$(x_n)_{n\in\Bbb N}$ and $(y_n)_{n\in\Bbb N}$ in $X$ so that
$(x_n)_n$ is bounded,  $\|x_n-y_n\|\geq\varepsilon_0>0$ 
for all $n\in\Bbb N$ and
$$
K_p(x_n,y_n)
=
2^{p-1}\|x_n\|^p + 2^{p-1}\|y_n\|^p - \|x_n+y_n\|^p
\leq
\frac {1}{n}\|x_n-y_n\|^p.
\eqno(\nfn)
$$
Moreover, without loss of generality we can suppose that
$\|y_n\|\geq\|x_n\|>0$ for all $n$. We then consider 
$0<\beta_n=\frac {\|x_n\|}{\|y_n\|}\leq 1$. From \pfn{0} it follows
$$
0\leq 2^{p-1}\big({\beta_n}^p +1\big) - \big(\beta_n +1\big)^p
\leq
\frac {1}{n}\big({\beta_n +1}\big)^p
@>>n\to\infty> 
0.
$$
Hence,
$$
2^{p-1} \frac{{\beta_n}^p +1}{(\beta_n+1)^p}
@>>n\to\infty>1.
\eqno(\nfn)
$$
\par
From \pfn{0} it follows that, 
$
\lim_{n\to\infty} \beta_n=\lim_{n\to\infty} \frac {\|x_n\|}{\|y_n\|}=1
$.
Then, since the sequence $(x_n)_n$ is bounded, so is
$(y_n)_n$ and therefore we have that
$\lim_{n\to\infty}\big(\|x_n\|-\|y_n\|\big)=0$. But using \pfn{-1} again,
we obtain that the bounded sequences $(x_n)_n$ and $(y_n)_n$ verify 
$$
\lim_{n\to\infty}\Big( \|x_n\|- \big\|\frac{x_n+y_n}{2}\big\|\Big)
=
\lim_{n\to\infty} \big(\|x_n\|-\|y_n\|\big)=0.
$$
Nonetheless, by hypothesis we have that
$\|x_n-y_n\|\geq\varepsilon_0>0$ for all $n$.
That is a 
contradiction with the uniform convexity of the norm $\norm$.\cqd
\par
Another important fact is that $K_p$ is uniformly separating when 
the modulus of convexity of the norm $\norm$ is of power type $p$. 
This is a consequence from results
of \cite{H}. Indeed, for any pair $x,\,y\in X$ we have
the following stronger inequality
$$
K_p(x,y)
=
{2^{p-1}\| x \|}^p + 2^{p-1}{\| y \|}^p - {\| x+y \|}^p 
\geq
C_{\|\cdot\|} {\| x - y \|}^p,
$$
for some $0<C_{\|\cdot\|}\leq 1$ (for instance, see \cite{C} Lemma 3.1).
\par
Hence, using {\bf Proposition 6} together with the inequality \pfn{-2}
we deduce (i), (iii), (iv) and (v) of {\bf Theorem 1}. It remains to prove
the assertion (ii), for which we use {\bf Proposition 8}. 
\par
More precisely, 
we observe that in the decomposition \pfn{-3} $d_K(\cdot, y)$ 
is a convex function
for every $y\in X$. Moreover, for $p>1$ is easy to verify that $c_{K_p}$ 
is strongly coercive; that is,
$$
\frac{c_{K_p}(x)}{\|x\|}
=
2^{p-1}\|x\|^{p-1}
@>>x\to\infty> +\infty.
$$
\par
Therefore, by {\bf Proposition 8}(iii) the
regularity of $\Delta_n^p{f}=\Delta_{K_p,n}f$ 
can be deduced from the regularity
of $c_{K_p}=2^{p-1}\norm^p$.
\par
Recall now that for any norm $\norm$ on $X$, 
the fact of being US (resp. with modulus of
smoothness of power type $1+\alpha$) is equivalent to 
$\norm\in\unifder(X)$ (resp. $\norm\in\holder{\alpha}(X)$). 
Therefore, 
$c_{K_p}=2^{p-1}\norm^p\in\bdunifder(X)$ 
(or $c_{K_p}\in\bdholder{\alpha}(X)$)
whenever the norm $\norm$ is US 
(or with modulus of smoothness of power type $1+\alpha$). 
\par
In the last case of (ii), for  
a norm $\norm$ with modulus of smoothness of power type
$1+\alpha$ (or equivalently $\norm\in\holder{\alpha}(X)$),
we can achieve a smoother
behaviour of the sequence $(\Delta_n^p{f})$ by choosing the proper 
value of $p$: 
$(\Delta_n^{1+\alpha}{f})_n\subset\holder{\alpha}(X)$. This is a corollary
of {\bf Proposition 8}(ii) and the next lemma.
\par
\proclaim{Lemma 1.3}
If $\norm\in\holder{\alpha}(X)$ then 
$\norm^{1+\alpha}\in\holder{\alpha}(X)$.
\endproclaim
\par
\demo{Proof of the Lemma 1.3}
This fact relies strongly in the convexity and homogeneity of a norm.
Since it is clear that $\norm^{1+\alpha}\in\bdholder{\alpha}(X)$,
let $C>0$ be the $\alpha$-H\"older continuity constant of the derivative
of the norm $\norm$ in $B_X$. We shall show that the condition \pfn{-16}
holds true for $\norm^{1+\alpha}$. Take any $x,\,y\in X$ and denote
by $\omega$ the maximum of $\|x\|$ and $\|y\|$. The lemma is proved
by the next computation.
$$
\eqalign{
\|x+y\|^{1+\alpha}+&\|x-y\|^{1+\alpha}-2\|x\|^{1+\alpha}
=
\cr
&\omega^{1+\alpha}
\Big(
\big\|\frac{x}{\omega}+\frac{y}{\omega}\big\|^{1+\alpha}
-
\big\|\frac{x}{\omega}\big\|^{1+\alpha}
+
\big\|\frac{x}{\omega}-\frac{y}{\omega}\big\|^{1+\alpha}
-
\big\|\frac{x}{\omega}\big\|^{1+\alpha}
\Big)
\leq
\cr
&\omega^{1+\alpha}2^{\alpha}C\big\|\frac{y}{\omega}\big\|^{1+\alpha}
=
2^{\alpha}C \|y\|^{1+\alpha}.
}
$$
\enddemo
\enddemo
By the above, this concludes the proof of {\bf Theorem 1}.
\cqd
\enddemo
\remark\nofrills{{\tensmc Acknowledgments.}}
\quad
The author wishes to thank Gilles Godefroy for 
his constant support and many fruitful conversations.
The author also wants to express his gratitude 
to the Department of Mathematics 
of the University of Missouri-Columbia,
where this work was developed. 
\endremark
%
%
%\vfill\break
%%%%%%%%%%%%%%%%%%%%%%%%%%%%%%%%%%%%%%%%%%%%%%%%%%%%%%%%%%%%%%%%%%%%%%%%%%%%%
%%%%%%%%% Bibliography
%
\refstyle{A}

\enddocument